\documentclass[11pt]{article}
\usepackage{amsthm}
\usepackage{ifthen,makeidx,amsmath,amssymb,epsf,
exscale,enumerate,array,xspace,fancybox,latexsym,isolatin1}
\usepackage[latin1]{inputenc}

\title{Sufficient Conditions for the Inversion Formula for the
  $k$-plane Radon
  Transform in $\Rn$.}
\author{Sine R. Jensen}
\date{}

\newcommand{\C}{{\mathbb C}}

\newcommand{\N}{{\mathbb N}}

\newcommand{\R}{{\mathbb R}}

\newcommand{\Rn}{{\R^n}}
\newcommand{\multp}{{\textbf p}}

\newcommand{\multl}{{\textbf l}}

\newcommand{\skK}{{\cal K}}
\newcommand{\skN}{{\cal N}}
\newcommand{\skO}{{\cal O}}
\newcommand{\skU}{{\cal U}}

\newcommand{\ra}{r^{\alpha}}
\newcommand{\xa}{x^{\alpha}_+}

\newcommand{\Ia}{I^{\alpha}}

\renewcommand{\Re}{{\operatorname {Re}}\,}
\renewcommand{\Im}{{\operatorname {Im}}\,}

\renewcommand{\epsilon}{\varepsilon}
\renewcommand{\phi}{\varphi}

\swapnumbers
\theoremstyle{plain}
\newtheorem{lemma}{Lemma}[section]
\newtheorem{proposition}[lemma]{Proposition}
\newtheorem{theorem}[lemma]{Theorem}
\newtheorem{corollary}[lemma]{Corollary}
\theoremstyle{definition}
\newtheorem{definition}[lemma]{Definition}
\theoremstyle{remark}
\newtheorem{remark}[lemma]{Remark}

\newcounter{romer}
\newcommand{\romer}{\begin{list}{(\Roman{romer})}{\usecounter{romer}}}
\newcommand{\slutromer}{\end{list}}

\begin{document} 
\def\appendixname{Appendix}
\maketitle

{\footnotesize\textsc{Abstract}: The inversion theorem
  (\ref{eq:formel}) for the
$k$-plane Radon transform in $\Rn$ is often stated for Schwartz functions,
cf. \cite[p.110]{GGA}, and lately for smooth functions on $\Rn$ fulfilling that $f(x)=O(|x|^{-N})$ for some $N>n$,
cf. \cite[Thm. I.6.2]{Helgason}. In this paper it will be shown,
that it suffices to require that $f$ is locally H\"older continuous
and $f(x)=O(|x|^{-N})$ for some $N>k$ ($N$ not necessarily an
integer) in order for (\ref{eq:formel}) to hold, and that the same
decay on $f$ but $f$ only continuous implies an inversion formula only
slightly weaker than (\ref{eq:formel}).}

{\small
\section*{Introduction}

An important area in the theory of the k-plane Radon transform on $\Rn$ is the
inversion theorems, which gives
explicit formulas by which one can recover a function from its
$k$-plane transform. Here we shall consider the formula
\begin{equation}
\label{eq:formel}
  f=(4\pi)^{-\frac{k}{2}}\frac{\Gamma\left(\frac{n-k}{2}\right)}
  {\Gamma\left(\frac{n}{2}\right)}I^{-k}(\hat{f})\,\check{},
\end{equation}
where ``$\quad\hat{}\quad$'' denotes the $k$-plane transform and ``$\quad\check{}\quad$'' the
dual transform, while $I^{-k}$ is a Riesz potential, cf. Section 4. It
will be shown in this paper, that the formula holds for all functions
in the space $C(k,n)$ (see Definition \ref{def:Holder3}.), and that
the formula with $I^k$ replaced by $\lim_{\alpha\to-k_+} I^\alpha$
holds if $f\in C_a(\R^n)$ for some $a>k$.

Notice, that the decay requirement of $C(k,n)$ ($f(x)=O(|x|^{-N})$ for
some $N>k$) on its member functions
is, in some sense, the weakest possible in order for an inversion
formula to hold: A sufficient condition for the integral in the
$k$-plane transform of a continuous function $f$ to be convergent is,
that for every $k$-plane there exists an $\epsilon>0$ such
that $f(x)=O(|x|^{-k-\epsilon})$ \textit{on this $k$-plane}. However
this non-uniform decrease of $f$ is not enough to make the inversion
formula valid. In \cite{Zalcman}, Zalcman shows the existence of a
smooth function $f\neq0$ on $\R^2$ satisfying $f(x)=O(|x|^{-2})$
\textit{on every line}, which nonetheless has $\hat{f}=0$. For further examples see e.g. \cite{A} and \cite{AG}.

The proof in this paper of the inversion formula is rooted in the
basic definition of the Riesz potential,
$I^\alpha$ ($\alpha\in\C$), which is
\begin{displaymath}
(I^\alpha f)(x)=\frac{1}{H_n(\alpha)}\int_\Rn f(y)|x-y|^{\alpha-n}\,dy.
\end{displaymath}
Here $H_n$ is a certain meromorphic function.
If $f$ is continuous and $O(|x|^{-a})$ for some $a>0$, the integral
converges if $0<\Re\alpha<a$. For values of $\alpha$ with
$\Re\alpha\leq0$, the Riesz potential can, depending on the regularity
of $f$, be defined by analytic
continuation (see e.g. \cite[sec. 10.2,
10.7]{Rubinbog} for various ways of performing this extension). The key to the proof of the inversion formula is the identity
$I^{-k}(I^kf)=f$, which will be established exactly for $f$ in $C(k,n)$.

Inversion formulas for the Radon transform of $L^p$-functions also
exists, but then the interpretation of the Riesz potentials is quite
different. Examples can be found e.g. in \cite{Rubin_e} where Rubin
verifies two inversion formulas for the case $k=n-1$. One of them is of
the same nature as (\ref{eq:formel}), and the other is of the type,
where a suitably interpreted Riesz potential in applied
\textit{before} the dual transform instead of after. The last
mentioned variant of inversion formula is in \cite{Rubin_b} proved for
$L^p$-functions in the case of a general $k$ under the assumption that
$1\leq p<\frac{n}{k}$. It is interesting to note, that given $f\in
C(\Rn)$ such that it is $O(|x|^{-N})$, then $f\in L^p(\Rn)$ 
when $-Np<-n$, i.e. $p>\frac{n}{N}$. Thus Rubin's inversion formula can
be used on this $f$ when there exists a $p\geq1$ with
$\frac{n}{N}<p<\frac{n}{k}$, e.i. when $k<N$ which is precisely the
decay condition in the inversion theorem of this paper.

The paper follows the lines of Helgason's exposition \cite[Chap.V
§5]{Helgason}: After the preliminaries, we study in Section 2 the analytic
continuation of the map $\alpha\mapsto\xa(f)=\frac{1}{\Gamma(\alpha+1)}\int_0^\infty
f(x)x^\alpha\,dx$. In Section 3 we use this to study the maps
$\alpha\mapsto\ra(f)=\frac{1}{\Gamma(\alpha+1)}\int_{\Rn}f(x)|x|^\alpha\,dx$,
and in Section 4 we introduce Riesz potentials and
establish the identity $I^{-k}(I^{k}f)=f$. Finally, in Section 5, we
prove the two versions of the inversion formula.

The inversion formula in (\ref{eq:formel}), expressed as it is in terms of
Riesz potentials, holds for $k$ both
odd and even. If $k$ is even it is well-known, that a similar inversion
formula can be established using the Laplacian instead of Riesz
potentials (see e.g. \cite[p. 29]{Helgason}. Section 6 contains a brief discussion of the possible impact of the
main result of the paper on the domain of this formula.
 
\section{Preliminaries}

For each $a>0$ and $n\in\N$ we make the following definitions:

\begin{definition}
Define the function space $C_a(\Rn)$ by
\begin{displaymath}
C_a(\Rn)=\{f\in C(\Rn)\mid f(x)=O(|x|^{-a})\}.
\end{displaymath}
\end{definition}

\begin{definition}
\label{def:Holder1}
For each $l\in\N_0=\N\cup\{0\}$, $0<\epsilon<1$ and $x\in\Rn$ define the space $C^{l+\langle\epsilon\rangle,x}(\Rn)$ as
the set of functions $f$ on $\Rn$ such that $f$ is $C^l$ in some neighborhood $\skO$ of $x$ with each 
$l$'th order derivative of $f$ H\"older continuous of index $\epsilon$ in that neighborhood, i.e.
\begin{equation}
\label{eq:Holder}
\exists M>0\forall
x_1,x_2\in\skO\forall\multl\in\N_0^n,|\multl|=l:|(\partial^\multl
f)(x_1)-(\partial^\multl f)(x_2)|\leq
M|x_1-x_2|^\epsilon.
\end{equation}
Put
\begin{itemize}
\item $C^{l+\langle\epsilon\rangle}(\Rn)
=\bigcap_{x\in\Rn}C^{l+\langle\epsilon\rangle,x}(\Rn)\subset C^l(\Rn)$,
\item $C^{l+}(\Rn)=\bigcap_{x\in\Rn}\bigcup_{\epsilon>0}C^{l+\langle\epsilon\rangle,x}\subset C^l(\Rn)$
\end{itemize}
and 
\begin{itemize}
\item $C^{l+\langle\epsilon\rangle,x}_a(\Rn)
=C^{l+\langle\epsilon\rangle,x}(\Rn)\cap C_a(\Rn)$,
\item
  $C^{l+\langle\epsilon\rangle}_a(\R^n)=C^{l+\langle\epsilon\rangle}(\R^n)\cap C_a(\Rn)$
\item $C^{l+}_a(\Rn)=C^{l+}(\Rn)\cap C_a(\R^n)$
\item $C^l_a(\Rn)=C^l(\Rn)\cap C_a(\Rn)$
\end{itemize}
\end{definition}

\begin{definition}
\label{def:Holder3}
Finally define for each $k\in\{1,\ldots,n-1\}$ the space $C(k,n)$ as the set of
functions $f$, such that $f\in C^{0+}_{k+\delta}(\Rn)$ for some
$\delta>0$. I.e. $f\in C(k,n)$ exactly when $f$ is
$O(|x|^{-k-\delta})$ for some $\delta>0$, and there for each $x\in\Rn$
exists a neighborhood $\skO$ and an $\epsilon$, $0<\epsilon<1$, such
that $|f(x_1)-f(x_2)|/|x_1-x_2|^{\epsilon}$ is bounded for $x_1,x_2\in\skO$.
\end{definition}

\vspace{0.1cm}

\begin{tabular}{|l|}
\hline
From now on, when the symbols $a$, $n$, $l$ and $\epsilon$ are used,
the assumption \\ will be $a>0$, $n\in\N$, $l\in\N_0$ and $0<\epsilon<1$,
unless otherwise mentioned.\\
\hline
\end{tabular}

\vspace{0.1cm}

\section{The map $\alpha\mapsto\xa(f)$}

\begin{definition}
For each $\alpha\in\C$ with $-1<\Re\alpha<a-1$ define the
map \hbox{$\xa:C_a(\R)\to\C$} by
\begin{equation}
  \label{eq:x_alpha}
  \xa(f)=\frac{1}{\Gamma(\alpha+1)}\int^\infty_0 f(x)x^\alpha dx.
\end{equation}
\end{definition}

\begin{remark}
\label{rem:x_alpha}
The map $\xa$ is well-defined since $-1<\Re\alpha$
and $f\in C(\R)$ makes the integrand integrable at 0, while
$\Re\alpha<a-1$ and $f(x)=O(|x|^{-a})$ makes it integrable at
$\infty$. Note, that the $\Gamma$-function is a non-vanishing
meromorphic function with poles in $-\N_0$ and 
\begin{equation}
\label{eq:Gamma}
\lim_{\alpha\to k}(\alpha-k)\Gamma(\alpha) =\frac{(-1)^{-k}}{(-k)!}, \qquad k\in-\N_0.
\end{equation}
\end{remark}

\begin{proposition}
\label{prop:x_alpha}
Let $f\in C^{l+\langle\epsilon\rangle,0}_a(\R)$. Then the map $\alpha\mapsto\xa(f)$, defined on 
\begin{displaymath}
  \{\alpha\in\C\mid-1<\Re\alpha<a-1\},
\end{displaymath}
can be (uniquely) extended to a holomorphic map on
\begin{displaymath}
  \{\alpha\in\C\mid-l-\epsilon-1<\Re\alpha<a-1\}.
\end{displaymath}
This map will likewise be denoted $\alpha\mapsto\xa(f)$. We have
\begin{equation}
\label{eq:residuum}
  \xa(f)=(-1)^{(-\alpha-1)}f^{(-\alpha-1)}(0),\qquad\textrm{when }\alpha\in\{-l-1,\ldots,-1\}.
\end{equation}
\end{proposition}

\begin{proof}
The integral in (\ref{eq:x_alpha}) is not necessarily convergent in
0, when $\alpha\leq-1$. But if we put
\begin{displaymath}
  A(x)=f(x)-\sum_{k=0}^l\frac{f^{(k)}(0)}{k!}x^k
\qquad \textrm{and}\qquad
  B(\alpha)=\sum_{k=0}^l\frac{f^{(k)}(0)\rho^{\alpha+k+1}}{k!(\alpha+k+1)}.
\end{displaymath}
then, by calculating the integrals, one
realizes that
\begin{equation}
\label{eq:x_alpha_m}
  \xa(f)=\frac{1}{\Gamma(\alpha+1)}(\int_0^\rho x^\alpha A(x)\,dx+\int_\rho^\infty x^\alpha f(x)\,dx+B(\alpha)),
\end{equation}
is an extension, cf. \cite[p.57]{GS}. Here $0<\rho<1$ fulfills
$\overline{B(0,\rho)}\subset\skO$, where $\skO$ is a neighborhood
of 0 in which $f^{(l)}$ is H\"older continuous. This extension is
well-defined on
\begin{displaymath}
  S=\{\alpha\in\C\setminus{-\N}\mid-l-\epsilon-1<\Re\alpha<a-1\}.
\end{displaymath}
To show this, only the first term needs
thought. Since $f\in C^{l}(\skO)$, there exists, according to Taylors
theorem, 
for any $x\in B(0,\rho)$ a $y$ between 0 and $x$, such that
\begin{equation}
\label{eq:taylor}
   f(x)=\sum_{k=0}^{l}\frac{f^{(k)}(0)}{k!}x^k+\frac{f^{(l)}(y)-f^{(l)}(0)}{l!}x^l.
\end{equation}
Because $f^{(l)}$ is H\"older continuous of index $\epsilon$ in $\skO$ we
therefore have
\begin{eqnarray}
  \label{eq:taylorb}
  \int_0^\rho |x^\alpha A(x)|\,dx
\leq const\int_0^\rho x^{\Re\alpha+l+\epsilon}\,dx<\infty
\end{eqnarray}
since $\Re\alpha+l+\epsilon>-l-\epsilon-1+l+\epsilon=-1$.

Let $\alpha_0\in S$ be given. To
  show that $\alpha\mapsto\xa(f)$ is holomorphic in
  $\alpha_0$, choose $\delta>0$ such that
  \begin{displaymath}
    B(\alpha_0,\delta)\subset\{\alpha\in\C\setminus{-\N}\mid-l-\epsilon-1+\delta<\Re\alpha<a-1-\delta\}.
  \end{displaymath}
Clearly $\alpha\mapsto B(\alpha)$ is holomorphic in $\alpha_0$. Thus
we only need to show, that the two integrals in (\ref{eq:x_alpha_m}) are holomorphic in
$\alpha_0$. This will follow from the theorems of Cauchy and Morera,
if it can be shown, that for any closed curve $\gamma$ in
$B(\alpha_0,\delta)$ the two integrals in each of the following
expressions can be interchanged:
\begin{displaymath}
  \int_\gamma\int_0^\rho x^\alpha A(x)
\qquad \textrm{and}\qquad
  \int_\gamma\int_\rho^\infty x^\alpha f(x).
\end{displaymath}
But for $x\in]0,\rho[$
$$
\sup_{\alpha\in B(\alpha_0,\delta)}|x^\alpha A(x)|\leq|A(x)|x^{-l-\epsilon-1+\delta},
$$
and this function is, as in (\ref{eq:taylorb}), integrable over
$]0,\rho[$.
For $x\in]\rho,\infty[$ we have the existence of a constant $c$
independent of $x$, such that
$$
\sup_{\alpha\in B(\alpha_0,\delta)}|x^\alpha f(x)|\leq c\, x^{a-1-\delta-a}=c\,x^{-1-\delta}.
$$

Now, let $m\in\{-l-1,\ldots,-1\}$ be given. Choose $\delta'>0$ such that
\begin{displaymath}
B(m,\delta')\setminus\{m\}\subset\{\alpha\in\C\setminus-\N\mid-l-\epsilon-1+\delta'<\Re\alpha<a-1-\delta'\}.
\end{displaymath}
As before we have for $\alpha\in B(m,\delta')$, that 
\begin{equation} 
  \label{eq:konstant}
  |\int_0^\rho x^\alpha A(x)\,dx|\leq C<\infty \qquad \textrm{and} \qquad |\int_\rho^\infty x^\alpha f(x)\,dx|\leq K<\infty,
\end{equation}
where the constants $C$ and $K$ are independent of $\alpha$. Thus for
$\alpha\to m$ we have
\begin{displaymath}
  (\alpha-m)\int_0^\rho x^\alpha
  A(x)\,dx\to0\qquad\textrm{and}\qquad(\alpha-m)\int_\rho^\infty x^\alpha
  f(x)\,dx\to0.
\end{displaymath}
Now (\ref{eq:residuum}) follows from (\ref{eq:x_alpha_m}) and (\ref{eq:Gamma}).
\end{proof}

\begin{remark}
\label{rem:no_Holder_x}
With the H\"older continuity condition on the derivatives of $f$
replaced by ordinary continuity, the inequality in (\ref{eq:taylorb}) changes to
$$
\int_0^\rho |x^\alpha A(x)|\,dx
\leq const\int_0^\rho x^{\Re\alpha+l}\,dx.
$$
Thus when $f\in C^l_a(\R)$, the extension of $\alpha\to\xa(f)$ still
exists but only on
$$
\{\alpha\in C|-l-1<\Re\alpha<a-1\}.
$$ 
\end{remark}

\section{The map $\alpha\mapsto\ra(f)$}

\begin{definition}
For each $\alpha\in\C$ with $-n<\Re\alpha<a-n$
define the map \hbox{$\ra:C_a(\Rn)\to\C$} by
\begin{equation}
  \label{eq:r_alpha}  
  \ra(f)=\frac{1}{\Gamma(\alpha+n)}\int_\Rn|x|^\alpha f(x)\,dx.
\end{equation}
\end{definition}
\begin{remark}
As in Remark \ref{rem:x_alpha} it is seen, that $r^\alpha$ is
well-defined.
\end{remark}

We will express $\ra$ by $\xa$. To this end we
introduce the mean value function:

\begin{definition}
For any $f\in C(\Rn)$ let $M_f:\R\to\C$ denote the mean value function
of $f$ around 0 defined by
\begin{equation}
\label{eq:middel}
  M_f(t)=\frac{1}{\Omega_n}\int_{S^{n-1}}f(t\omega)\,d\omega.
\end{equation}
\end{definition}

\begin{remark}
\label{rem:middel}
Notice, that $t\mapsto M_f(t)$ is even, and that $M_f(0)=f(0)$.
\end{remark}

\begin{lemma}
\label{lemma:middel}
When $f$ is in $C^{l+\langle\epsilon\rangle,0}_a(\Rn)$ then $M_f$ is in $C^{l+\langle\epsilon\rangle,0}_a(\R)$.
\end{lemma}

\begin{proof}
Standard arguments.
\end{proof}

\begin{remark}
\label{rem:r_alpha}
Transition to polar coordinates in the defining expression
(\ref{eq:r_alpha}) for $\ra$ now gives $\ra(f)$ in terms of $\xa$: 
\begin{equation}
  \label{eq:ra_ved_xa}
  \ra(f)=\Omega_nx_+^{\alpha+n-1}(M_f),
\end{equation}
when $-1<\Re\alpha+n-1<a-1$, i.e. $-n<\Re\alpha<a-n$.
\end{remark}

\begin{proposition}
\label{prop:r_alpha}
Let $f\in C^{l+\langle\epsilon\rangle,0}_a(\Rn)$. Then the map $\alpha\mapsto\ra(f)$, defined on
\begin{displaymath}
  \{\alpha\in\C\mid-n<\Re\alpha<a-n\},
\end{displaymath}
can be (uniquely) extended to a holomorphic map on
\begin{displaymath}
  A=\{\alpha\in\C\mid-l-\epsilon-n<\Re\alpha<a-n\}.
\end{displaymath}
This map will likewise be denoted $\alpha\to\ra(f)$, and it satisfies
(\ref{eq:ra_ved_xa}). In specific
\begin{equation}
\label{eq:r_graenser}
  \ra(f)=\Omega_n (-1)^{-\alpha-n} M_f^{(-\alpha-n)}(0),\qquad\textrm{when }\alpha\in\{-l-n,\ldots,-n\}.
\end{equation}
\end{proposition}

\begin{proof}
Use (\ref{eq:ra_ved_xa}) as definition and apply
Proposition \ref{prop:x_alpha} using Lemma \ref{lemma:middel}.
\end{proof}

\section{Riesz Potentials}

\begin{definition}
The meromorphic function $H_n$ on $\C$ is defined by
\begin{displaymath}
  H_n(\alpha)=2^\alpha\pi^{\frac{n}{2}}\frac{\Gamma(\frac{\alpha}{2})}{\Gamma(\frac{n-\alpha}{2})}.
\end{displaymath}
\end{definition}

\begin{remark}
Note that $H_n$ has simple poles at each $\alpha\in-2\N_0$ and a zero
in each $\alpha\in n+2\N_0$.
\end{remark}
\begin{definition}
We put $\C_n=\C\setminus(n+2\N_0)$.
\end{definition}
\begin{definition} 
\label{def:Riesz}
For each $x\in\Rn$, $f\in C_a(\Rn)$, and $\alpha\in\C_n$ with
$0<\Re\alpha<a$ the 
{\em {$\alpha$th Riesz
  potential, $I^\alpha$, of $f$ at $x$}} is defined as
\begin{equation}
  \label{eq:defRiesz}
  (I^\alpha f)(x)=\frac{1}{H_n(\alpha)}\int_\Rn f(y)|x-y|^{\alpha-n}\,dy
=\frac{1}{H_n(\alpha)}\int_\Rn f(x-y)|y|^{\alpha-n}\,dy.
\end{equation}
\end{definition} 

\begin{remark} \label{rem:welldef}
As in Remark \ref{rem:x_alpha} it is seen, that $I^\alpha f(x) $ is
well-defined. Comparing with the defining expression (\ref{eq:r_alpha}) for $\ra$ we
see, that
\begin{equation}
  \label{eq:overensstemmelse}
  (I^\alpha
  f)(x)=\frac{\Gamma(\alpha)}{H_n(\alpha)}r^{\alpha-n}(\tau_xf)
\end{equation}
where $\tau_x f(y)=f(x-y)$.
\end{remark}

\begin{proposition}
\label{prop:I_alpha}
Let $x\in\Rn$ be given. Assume that $f\in C^{l+\langle\epsilon\rangle,x}_{a}(\Rn)$. Then the map
$\alpha\mapsto (I^\alpha f)(x)$, defined on the set
\begin{displaymath}
  \{\alpha\in\C_n\mid0<\Re\alpha<a\},
\end{displaymath}
can be (uniquely) extended to a meromorphic map on
\begin{displaymath}
  B=\{\alpha\in\C\mid-l-\epsilon<\Re\alpha<a\}.
\end{displaymath}
This map will likewise be denoted $\alpha\to(I^\alpha f)(x)$. It
satisfies (\ref{eq:overensstemmelse}) for\\ \noindent \hbox{$\alpha\in
B\setminus((-\N_0)\cup(n+2\N_0))$}. The poles, which are all simple,
are in
$$
(n+2\N_0)\cup B.
$$
\end{proposition}

\begin{proof}
Use (\ref{eq:overensstemmelse}) as definition and apply
Proposition \ref{prop:r_alpha} to obtain a (unique) meromorphic extension to
$\{\alpha\in\C\mid-l-\epsilon<\Re\alpha<a\}$. The possible poles are
those of
$\frac{\Gamma(\alpha)}{H_n(\alpha)}=\frac{1}{2}\pi^{-\frac{n+1}{2}}\Gamma(\frac{n-\alpha}{2})\Gamma(\frac{\alpha+1}{2})$.
They are $\alpha\in2\N_0+n$ and $\alpha\in-2\N_0-1$, all simple. When
$\alpha\in(-2\N_0-1)\cap B$
it follows from (\ref{eq:r_graenser}), that
$$
r^{\alpha-n}(f)=\Omega_n (-1)^{-\alpha}M_f^{(-\alpha)}(0)=0,
$$
since $M_f$ in an even function. Thus $\alpha$ is a removable singularity.
\end{proof}


\begin{lemma} 
\label{lemma:If}
Let $f\in C^{0+}_{a}$. Then $x\mapsto(I^0f)(x)$ is defined on all
of $\Rn$ and 
$$
I^0f=f.
$$
\end{lemma}
\begin{proof}
It follows from Proposition \ref{prop:I_alpha}, that $x\mapsto(I^0f)(x)$ is
defined on all of $\Rn$. Since
$$
\lim_{\alpha\to 0}\alpha H_n(\alpha)=\frac{2\pi^{\frac{n}{2}}}{\Gamma(\frac{n}{2})}=\Omega_n,
$$
it follows from Proposition \ref{prop:I_alpha}, (\ref{eq:r_graenser}),
(\ref{eq:Gamma}) and Remark \ref{rem:middel}, that

\begin{equation}
\label{eq:If}
  (I^0f)(x)=\lim_{\alpha\to0}\frac{\alpha\Gamma(\alpha)
   }{\alpha H_n(\alpha)} r^{\alpha-n}(\tau_xf)
=M_{\tau_xf}(0)=f(x).
\end{equation}
\end{proof}

\begin{lemma}
\label{lemma:decay}
Let $f\in C_a(\Rn)$. Let $\alpha\in\C$ with
$0<\Re\alpha<\min(a,n)$ be given. Then  
$$
\Ia f\in C_{b-\Re\alpha}(\Rn),
$$
for any $b$ with
$\Re\alpha<b\leq\min(a,n)$ if $a\neq n$, and for any $b$ with
$\Re\alpha<b<n$ if $a=n$.
\end{lemma}

\begin{proof}
See \cite[Prop. V.5.8.]{Helgason} with natural modifications to
the proof in case $a=n$.
\end{proof}

\begin{proposition} \label{prop:sammensaetning}
Let $f\in C_a(\Rn)$. For any pair $\alpha,\beta\in\C$ satisfying
\begin{displaymath}
\label{eq:restriktiv}
  \Re\alpha>0\quad\textrm{and}\quad\Re\beta>0\quad\textrm{and}\quad\Re(\alpha+\beta)<\min(a,n)
\end{displaymath}
we have
\begin{equation}
\label{eq:modulo_en_konstant}
  I^\alpha I^\beta f=I^{\alpha+\beta}f.
\end{equation}
\end{proposition}

\begin{remark}
Refer to e.g. \cite[p. 43ff]{Landkof} or \cite[Satz 9]{Ortner} in order to
see how, when dealing with Riesz potentials as distributions,
(\ref{eq:modulo_en_konstant}) can be expressed as a convolution of
distributions. The distribution approach can prove Proposition
\ref{prop:sammensaetning} for a smaller class of functions.
\end{remark}

\begin{proof}
That $0<\Re\beta<\min(a,n)$ implies two things. First we get from Remark \ref{rem:welldef}, that $I^\beta f$ is well-defined and given by 
\begin{displaymath}
(I^{\beta}f)(z)=\frac{1}{H_n(\beta)}\int_\Rn f(y)|z-y|^{\beta-n}\,dy.
\end{displaymath}
Secondly, we get the usage of Lemma \ref{lemma:decay} from which follows, that 
$$
I^\beta f\in C_{b-\Re\beta}(\Rn),
$$
where $b$ is chosen such that $\Re(\alpha+\beta)<b<\min(a,n)$.
Thus, because $0<\Re\alpha<b-\Re\beta$, $I^{\alpha}(I^{\beta}f)$ is well-defined and given by
\begin{eqnarray}\label{eq:dobbeltintegral}
I^{\alpha}(I^{\beta}f)(x)
&=&\frac{1}{H_n(\alpha)}\int_\Rn (I^\beta f)(z)|x-z|^{\alpha-n}\,dz\nonumber\\
&=&\frac{1}{H_n(\alpha)}\frac{1}{H_n(\beta)}\int_\Rn\int_\Rn f(y)|z-y|^{\beta-n}\,dy|x-z|^{\alpha-n}\,dz.
\end{eqnarray}
To show, that the order of integration can be interchanged, consider
the expression
\begin{equation}
  \label{eq:sammensat}
  \int_\Rn |f(y)|\int_\Rn |z-y|^{\Re\beta-n}|x-z|^{\Re\alpha-n}\,dz\,dy.
\end{equation}
By substituting $v=\frac{x-z}{|x-y|}$ in the inner integral and using
the rotation invariance of the Lebesgue measure, this expression is
rewritten as
\begin{displaymath}
  \int_\Rn |f(y)||x-y|^{\Re\alpha+\Re\beta-n}\,dy\int_\Rn|e-v|^{\Re\beta-n}|v|^{\Re\alpha-n}\,dv,
\end{displaymath}
where $e$ is an arbitrary fixed unit vector. Now
$0<\Re(\alpha+\beta)<a$ makes the $y$-integral convergent. That the
$v$-integral is convergent can be seen easily. Finally, it can be
shown, e.g. using Fourier transform as in \cite[p. 117-118]{Stein}, that 
$$
\int_\Rn|e-v|^{\beta-n}|v|^{\alpha-n}\,dv=\frac{H_n(\alpha)H_n(\beta)}{H_n(\alpha+\beta)}.
$$
\end{proof}

\begin{remark}
\label{rem:decomposition}
Let $x_0\in\Rn$ be given. In what follows, we will often decompose
a given function $f$ on $\Rn$ as $f=f_1+f_2$, where $f_1=(1-\chi)f$
and $f_2=\chi f$ for some compactly supported $C^\infty$-function
$\chi$ with $\chi(x)=1$ in some neighborhood of $x_0$. Note that $f_1$ and
$f_2$ have the same regularity as $f$, but $f_1$ is 0 in the
neighborhood of $x_0$ and $f_2$ has compact support.
\end{remark}

\begin{lemma}
\label{lemma:decomposition}
Let $f\in C^l_a(\Rn)$. Let $\alpha\in \C_n$ with $0<\Re\alpha<a$ and
$x_0\in\Rn$ be given. Write $f=f_1+f_2$ as in Remark
\ref{rem:decomposition}. Then $I^\alpha f_1$ is smooth at $x_0$, and
$\Ia f_2\in C^l(\Rn)$ with
$\partial^\multp(\Ia f_2)=\Ia(\partial^\multp f_2)$ for any
$\multp\in\N_0^n$ with $|\multp|\leq l$.
\end{lemma}

\begin{proof}
Assume $\multp\in\N_0^n$ to be given. Choose $\delta>0$ such that
$f_1=1$ in
$B(x_0,2\delta)\subset\skU$. Then for any $x\in B(x_0,\delta)$
\begin{eqnarray*}
|f_1(y)\partial_x^\multp|x-y|^{\alpha-n}|&\leq&c|f_1(y)||x-y|^{\Re\alpha-n-|\multp|}\\
&\leq&c'1_{\Rn\setminus B(x_0,2\delta)}(y)(|y|+1)^{-a}|y-x_0|^{\Re\alpha-n-|\multp|},
\end{eqnarray*}
since $\frac{1}{2}|x_0-y|\leq|x-y|\leq2|x_0-y|$ for
$y\notin B(x_0,2\delta)$. Here $c'$ does not depend on $x$. Since $-a+\Re\alpha-n-|\multp|<-a+a-n=-n$,
this gives us an integrable majorant of
\hbox{$\partial_x^\multp(f_1(y)|x-y|^{\alpha-n})$} and it is independent of $x$. 

To deal with $I^\alpha f_2$ assume $|\multp|\leq l$ and let $\skN$ be any
bounded subset of $\Rn$. Let $x\in\skN$. Then
\begin{eqnarray}
\label{eq:diffint}
|\partial_x^\multp (f_2(x-y)|y|^{\alpha-n})|&=&|(\partial^\multp
f_2)(x-y)||y|^{\Re\alpha-n}|\nonumber\\
&\leq&
\sup|\partial^\multp f_2|\,1_{\skN+(-\skK)}(y)|y|^{\Re\alpha-n},\qquad\forall y\in\Rn.
\end{eqnarray}
Since $\Re\alpha>0$ this is an
integrable majorant of \hbox{$\partial_x^\multp
  (f_2(x-y)|y|^{\alpha-n})$} and it is 
independent of $x$. Thus $\partial^\multp (I^\alpha f_2)$ exists in
$\skN$, $\skN$ arbitrary, and thus in all of $\Rn$, and we see from
(\ref{eq:diffint}) that
$\partial^\multp(\Ia f_2)=\Ia(\partial^\multp f_2)$.
\end{proof}

\begin{lemma} \label{lemma:regularity1}
Let $f\in C^l_a(\Rn)$. Let $\alpha\in\C_n$ with $0<\Re\alpha<a$ be
given. Then
\begin{equation}
\label{eq:c}
I^\alpha f\in C^l(\Rn)
\end{equation}
and for any $x\in\Rn$ and $0<\epsilon<1$
\begin{equation}
\label{eq:Hc}
f\in C^{l+\langle\epsilon\rangle,x}(\Rn)\Rightarrow I^\alpha f\in C^{l+\langle\epsilon\rangle,x}(\Rn).
\end{equation}
\end{lemma}

\begin{proof}
Let $x_0\in\Rn$ be given. Write $f=f_1+f_2$ as in Remark
\ref{rem:decomposition}. From the preceeding lemma $\Ia f_1$ is smooth
at $x_0$ and $\Ia f_2\in C^l(\Rn)$. Thus (\ref{eq:c}) holds. 

Assume now, that $f\in C^{l+\langle\epsilon\rangle,x_0}(\Rn)$. Let
$\multl\in\N_0^n$ with $|\multl|=l$ be given. To show the H\"older
continuity of $\partial^\multl (I^\alpha f_2)$ ($=\Ia(\partial^\multl
f_2)$ according to Lemma \ref{lemma:decomposition}), let $\skK$ be a
compact neighborhood of $x_0$ in which $\partial^multl f$ is H\"older
continuous of index $\epsilon$ and assume $\chi$ in the decomposition
$f=f_1+f_2=(1-\chi)f+\chi f$ to have $\skK$ as its support. Then
$\partial^\multl f_2$ is H\"older continuous of index $\epsilon$ in all of
$\Rn$, so for any bounded neighborhood $\skN$ of $x_0$ and any
$x_1,x_2\in\skN$ 
\begin{eqnarray*}
&&|\partial^\multl (I^\alpha f_2)(x_1)-\partial^\multl (I^\alpha f_2)(x_2)|
\\
&\leq&\frac{1}{H_n(\alpha)}\int_{\skN+(-\skK)}|\partial^\multl
f_2(x_1-y)-\partial^\multl f_2(x_2-y)||y|^{\Re\alpha-n}\,dy
\leq M'|x_1-x_2|^\epsilon
\end{eqnarray*}
for some $M'>0$.
\end{proof}

\begin{lemma}
\label{lemma:regularity2}
Let $f\in C_a(\Rn)$. Let $\alpha\in\C_n$ with $\Re\alpha=1$ be
given. If $a>1$ then 
$$
f\in C^{l+\langle\epsilon\rangle,x}(\Rn)\Rightarrow \forall\epsilon',0<\epsilon'<\epsilon:I^\alpha f\in C^{(l+1)+\langle\epsilon'\rangle,x}(\Rn)
$$
for any $x\in\Rn$ and $0<\epsilon<1$.
\end{lemma}
\begin{proof}
Let $x_0\in\Rn$ be given. Write $f=f_1+f_2$ as in Remark
\ref{rem:decomposition}. Then from Lemma \ref{lemma:decomposition} $\Ia
f_1$ is smooth at $x_0$, so only $\Ia f_2$ needs thought.

Pick $\multp\in\N^n_0$ with
$|\multp|=l+1$. Write $\multp=\multl+e_i$ for some $\multl\in\N^n_0$
with $|\multl|=l$, and some
$e_i=(0,\ldots,0,1,0,\ldots,0)$. Let $\skK$ be a compact neighborhood
of $x_0$ in which $\partial^l f$ is H\"older continuous and assume
$\chi$ in the decomposition $f=f_1+f_2=(1-\chi)f+\chi f$ to have
$\skK$ as its support. Put $g=\partial^\multl f_2$. Then $g$ 
is H\"older continuous of index $\epsilon$ in all of $\Rn$ and has
support in $\skK$. What needs to be shown is, that $\partial^\multp
(I^\alpha f_2)=\partial_i\partial^\multl 
(I^\alpha f_2)=\partial_i (I^\alpha g)$ (Lemma \ref{lemma:decomposition})
exists and is H\"older
continuous in a neighborhood of $x_0$.

Let $B$ be a symmetric, bounded neighborhood of 0 such that
$\skK\subset B+x=B_x$ for all $x$ in some bounded, open neighborhood
$\skO$ of $x_0$. Let $\beta\in\C$ with $1<\Re\beta<2$ be given. Then for any $x\in\skO$
\begin{eqnarray}
\partial_iI^\beta g(x)
&=&c_n(\beta)\int_{B_x}g(y)(x_i-y_i)|x-y|^{\beta-n-2}\,dy\nonumber\\
&=&c_n(\beta)\int_{B}g(x-y)y_i|y|^{\beta-n-2}\,dy
\end{eqnarray}
where $c_n(\beta)=\frac{\beta-n}{H_n(\beta)}$ and where the integral exists since $\Re\beta>1$ and $B$ is bounded. Furthermore, using the H\"older continuity of $g$
$$
\int_B|(g(x-y)-g(x))y_i|y|^{\alpha-n-2}|\,dy\leq M\int_B|y|^{-n+\epsilon}\,dy<\infty,
$$
i.e. the integral $\int_B(g(x-y)-g(x))y_i|y|^{\alpha-n-2}\,dy$
exists. Using the symmetry of $B$ we get
\begin{eqnarray}
\label{eq:limit}
&&|\frac{1}{c_n(\beta)}\partial_iI^\beta g(x)-\int_B(g(x-y)-g(x))y_i|y|^{\alpha-n-2}\,dy|\nonumber\\
&\leq&\int_B|(g(x-y)-g(x))y_i(|y|^{\beta-n-2}-|y|^{\alpha-n-2})|\,dy+|g(x)\int_B
y_i|y|^{\beta-n-2}\,dy|\nonumber\\
&\leq&c'\int_B||y|^{\beta-n-1+\epsilon}-|y|^{\alpha-n-1+\epsilon}|\,dy
\end{eqnarray}
for some $c'>0$. Now notice that when $n=1$, then $c_n$ has a removable
singularity at $\beta=1$, so that for any value of $n\in\N$, $c_n$ is
bounded and bounded away from 0 in a small enough neighborhood of
$\alpha$, i.e. $\lim_{\beta\to\alpha} \frac{1}{c_n(\beta)}$ exists
and is not 0. Thus (\ref{eq:limit}) shows that in the limit where $\Re\beta>1$
$$
\lim_{\beta \to \alpha} \partial_i I^\beta g(x)=c_n(\alpha)\int_B(g(x-y)-g(x))y_i|y|^{\alpha-n-2}\,dy
$$ 
uniformly on $\skO$. So $\partial_i I^{\alpha}g$ does exist and
$$
\partial_i I^{\alpha}g=c_n(\alpha)\int_B(g(x-y)-g(x))y_i|y|^{\alpha-n-2}\,dy
$$
in all of $\skO$. Given $0<\epsilon'<\epsilon$ put
$s=\frac{\epsilon'}{\epsilon}$ and $t=1-s$. We then have for any $x_1,x_2\in\skO$, that
\begin{eqnarray*}
&&|\partial_i I^{\alpha}g(x_1)-\partial_i I^{\alpha}g(x_2)|\\
&=&c_n(\alpha)|\int_B(g(x_1-y)-g(x_1)-(g(x_2-y)-g(x_2)))y_i|y|^{\alpha-n-2}\,dy\\
&\leq&c_n(\alpha)\int_B|(g(x_1-y)-g(x_1))-(g(x_2-y)-g(x_2))|^t\\
&&\qquad|(g(x_1-y)-g(x_2-y))-(g(x_1)-g(x_2))|^s|y|^{-n}\,dy\\
&\leq&c_n(\alpha)\int_B(2M|y|^{\epsilon})^t(2M|x_1-x_2|^{\epsilon})^s|y|^{-n}\,dy\\
&=&c_n(\alpha)|x_1-x_2|^{\epsilon s}2M\int_B|y|^{\epsilon
  t-n}\,dy=M'|x_1-x_2|^{\epsilon'}
\end{eqnarray*}
for some constant $M'>0$.
\end{proof}

\begin{corollary} \label{cor:I_decrease}
Let $f\in C_{a}(\Rn)$. Let $\alpha\in\C$ with
$0<\Re\alpha<\min(a,n)$ be given. Then 
$$
f\in C^{l+\langle\epsilon\rangle,x}(\Rn)\Rightarrow \forall\epsilon',0<\epsilon'<\epsilon:I^\alpha f\in C^{(l+[\Re\alpha])+\langle\epsilon'\rangle,x}(\Rn)
$$
for any $x\in\Rn$ and $0<\epsilon<1$.
Here $[\Re\alpha]$ denotes the integer part of $\Re\alpha$.
\end{corollary}

\begin{proof}
Write $\alpha=\beta+[\Re\alpha]$. Then $0\leq\Re \beta
<1$. From Proposition \ref{prop:sammensaetning} combined with
Lemma \ref{lemma:decay}
$$
I^\alpha f=I^\beta(I^1(I^1(\ldots(I^1 f)\ldots))),
$$
$I^1$ applied $[\Re\alpha]$ times. The claim now follows from Lemma
\ref{lemma:regularity2} and Lemma \ref{lemma:regularity1}.
\end{proof}

\begin{proposition} \label{prop:sammensaetning_-k_k}
Let $k\in\{1,\ldots,n-1\}$ and $f\in C(k,n)$.
Then 
$$
I^{-k}(I^{k}f)=f.
$$
\end{proposition}
\begin{proof}
Let $x\in\Rn$ be given and choose $\delta$, $0<\delta<1$, such that
$f\in C_{k+\delta}(\Rn)$. From Proposition \ref{prop:I_alpha} it follows, that there exists an $\delta'$, $0<\delta'<1$, such that the map
$$
\alpha\mapsto(I^{\alpha+k}f)(x)
$$
is holomorphic in
$\{\alpha\in\C\mid-k-\delta'<\Re\alpha<\delta\}$. Since Lemma
\ref{lemma:decay} and Corollary \ref{cor:I_decrease} with $a=b=k+\delta$ ensures, that
\begin{displaymath}
  I^kf\in C^{k+}_{\delta}(\Rn),
\end{displaymath}
we likewise get from Proposition \ref{prop:I_alpha}, that there exists a $\delta''$, $0<\delta''<1$, such that the map
$$
\alpha\mapsto(I^{\alpha}(I^{k}f))(x)
$$
is well-defined and holomorphic in $\{\alpha\in\C\mid-k-\delta''<\Re\alpha<\delta\}$. Proposition \ref{prop:sammensaetning} gives us, that
$$
I^\alpha I^k f(x)=I^{\alpha+k}f(x),
$$
when $\alpha\in\{\alpha\in\C\mid 0 < \Re\alpha < \delta \}$. By
analytic continuation this identity then holds on all of
$\{\alpha\in\C\mid-k-\min(\delta',\delta'')<\Re\alpha<\delta\}$.
In particular, using Lemma \ref{lemma:If} with $a=k+\delta$
$$
I^{-k} I^k f(x)=I^0f(x)=f(x).
$$
\end{proof}

\section{The Inversion Formula for the Radon Transform}

Let $k\in\{1,\ldots,n-1\}$ be given. Let $f\in C_a(\Rn)$ for some
$a>k$. For the $k$-plane transform one arrives, by calculating, at
\begin{equation} \label{eq:udregning}
  (\hat{f})\,\check{}\,(x)=(4\pi)^{\frac{k}{2}}\frac{\Gamma\left(\frac{n}{2}\right)}{\Gamma\left(\frac{n-k}{2}\right)}(I^kf)(x),
\end{equation}
cf. \cite{Fuglede} or \cite[p.29]{Helgason}. This will be used in what
follows.

\begin{theorem}
  \label{thm:omvendingssaetningen}
Let $k\in\{1,\ldots,n-1\}$.
Assume, that 
$f\in C(k,n)$.
Then $f$ can be recovered from its $k$-plane transform by
\begin{displaymath}
  f=(4\pi)^{-\frac{k}{2}}\frac{\Gamma\left(\frac{n-k}{2}\right)}
  {\Gamma\left(\frac{n}{2}\right)}I^{-k}(\hat{f})\,\check{}.
\end{displaymath}
\end{theorem}

\begin{proof}
The claim follows from (\ref{eq:udregning}) by means of Proposition \ref{prop:sammensaetning_-k_k}.
\end{proof}

\begin{remark}
Any differentiable function will also be locally H\"older continuous
(but the inverse implication is not true). Thus the H\"older condition
could in the entire paper have been replaced by demanding all
functions to be one more time continuously differentiable. E.g. Theorem
\ref{thm:omvendingssaetningen} is true for all $f\in C^1(\Rn)$
with $f(x) = O(|x|^{-k-\delta})$ for some $\delta>0$.
\end{remark}

An even lower regularity requirement on $f$ can be bought at a small price:

\begin{theorem}
\label{thm:alternativt}
Let $k\in\{1,\ldots,n-1\}$. Assume, that $f\in C_{k+\delta}$ for some
$\delta>0$. Then $f$ can be recovered from its $k$-plane transform by
\begin{displaymath}
  f=(4\pi)^{-\frac{k}{2}}\frac{\Gamma\left(\frac{n-k}{2}\right)}
  {\Gamma\left(\frac{n}{2}\right)}\lim_{s\to-k_+}
  I^{s}(\hat{f})\,\check{}.
\end{displaymath}
\end{theorem}

We will need the following lemma pointed out to me by Boris Rubin
(cf. \cite[Thm. I.2.6]{Samko}): 

\begin{lemma}
\label{lemma:alternativt}
  Let $f\in C_a(\R)$. Then
\begin{displaymath} 
  \lim_{s\to-1_+}x^s_+(f)=f(0).
\end{displaymath}
\end{lemma}

\begin{proof}
Let $\epsilon>0$ be given and choose $\delta$,
$0<\delta<1$, such that $|f(x)-f(0)|\leq\epsilon$ when
$|x|\leq\delta$. Write
\begin{displaymath}
x^s_+(f)=\frac{1}{\Gamma(s+1)}\left[
\int_0^\delta(f(x)-f(0))x^s\,dx+
\int_\delta^\infty f(x)x^s\,dx+
\int_0^\delta f(0)x^s\,dx\right].
\end{displaymath}
Since (when $ s>-1$),
\begin{displaymath}
| \frac{1}{\Gamma(s+1)} 
\int_0^\delta(f(x)-f(0))x^s\,dx|
\leq\frac{\epsilon}{\Gamma( s+2)}\delta^{ s+1}
\end{displaymath}
and (when $ s-a<-1$)
\begin{displaymath}
|\frac{1}{\Gamma(s+1)}\int_\delta^\infty
  f(x)x^s\,dx|
\leq\frac{c}{\Gamma( s+1)| s-a+1|}\delta^{ s-a+1}
\end{displaymath}
for some constant $c>0$ and
\begin{displaymath}
|\frac{1}{\Gamma(s+1)}\int_0^\delta
  f(0)x^s\,dx-f(0)|
\leq|
f(0)(\frac{\delta^{s+1}}{\Gamma(s+2)}-1
)|,
\end{displaymath}
$|x^s_+(f)-f(0)|$ can be estimated by e.g. some multiple of
$\epsilon$ when $s$ is sufficiently close to $-1$.
\end{proof}

Also a parallel of Corollary \ref{cor:I_decrease} and thus of
Lemma \ref{lemma:regularity2} for functions with the H\"older
continuity of the derivatives replaced by ordinary continuity is
needed:
\begin{lemma}
Let $f\in C^l_a(\Rn)$. Let $\alpha\in\C_n$ with $Re\alpha=1$ be
given. If $a>1$, then 
$$
I^\alpha f\in C^{l+\langle\epsilon\rangle}(\Rn)
$$
for any $0<\epsilon<1$.
\end{lemma}

\begin{proof}
Let $x\in\Rn$ and $0<\epsilon<1$ be given. Decompose
$f=f_1+f_2$ as in Remark \ref{rem:decomposition}. Then $\Ia f_1$
is smooth at $x$
according to Lemma \ref{lemma:decomposition}, so only $\Ia f_2$ needs
thought. 

From Lemma
\ref{lemma:decomposition} $\Ia f_2$ is in $C^l(\Rn)$ 
with $\partial^\multl\Ia f_2=\Ia\partial^\multl f_2$ for any
$\multl\in\N_0^n$ with $|\multl|=l$. The claim
is, that these derivatives are H\"older continuous of index
$\epsilon$ at $x$. Therfore pick $\skO$, a bounded neighborhood of $x$,
and $x_1,x_2\in\skO$. Since $f_2$ has compact support $\skK$, there
exists $c>0$ such that
$$
|\partial^\multl\Ia f_2(x_1)-\partial^\multl\Ia
f_2(x_2)|\leq c\int_\skK\big||x_1-y|^{\alpha-n}-|x_2-y|^{\alpha-n}\big|\,dy.
$$
Thus it suffices to prove the existence of a constant $C>0$ (independent of
$x_1$ and $x_2$) such that
\begin{equation}
\label{eq:integralK}
\int_\skK\big||x_1-y|^{\alpha-n}-|x_2-y|^{\alpha-n}\big|\,dy\leq C|x_1-x_2|^\epsilon.
\end{equation}
Put
$$
B_1=B(x_1,\frac{2}{3}|x_1-x_2|)\qquad
B_2=B(x_2,\frac{2}{3}|x_1-x_2|)\qquad
A=\skK\setminus(B_1\cup B_2).
$$
Then $\skK\subset B_1\cup B_2\cup A$, so (\ref{eq:integralK}) holds if
it can be proved with $\skK$ replaced by each of the three sets
$B_1$, $B_2$ and $A$. But since $|x_2-y|>\frac{1}{3}|x_1-x_2|$ when
$y\in B_1$
\begin{eqnarray*}
&&\int_{B_1}\big||x_1-y|^{\alpha-n}-|x_2-y|^{\alpha-n}\big|\,dy\\
&\leq&\int_{B_1}|x_1-y|^{1-n}\,dy+\int_{B_1}|x_2-y|^{1-n}\,dy\\
&\leq&\int_{B(0,\frac{2}{3}|x_1-x_2|)}
|y|^{1-n}\,dy+\int_{B(0,\frac{2}{3}|x_1-x_2|)}
(\frac{1}{3}|x_1-x_2|)^{1-n}\,dy \leq C_1|x_1-x_2|.
\end{eqnarray*}
An equivalent calculation can be done for the integral on $B_2$. Thus we
turn to the integral on $A$.

First let $y\in A$ with $|x_1-y|\neq|x_2-y|$ be given. Apply the mean
value theorem to the function $t\mapsto\Re t^{\alpha-n}$
on the interval with endpoints $|x_1-y|$ and $|x_2-y|$ to obtain the
existence of an $s_1\in]0,1[$ such that
$$
\big|\Re|x_1-y|^{\alpha-n}-\Re|x_2-y|^{\alpha-n}\big|\leq c'(s_1|x_1-y|+(1-s_1)|x_2-y|)^{-n}\big||x_1-y|-|x_2-y|\big|.
$$
Then apply the mean value theorem to the function $t\mapsto\Im
t^{\alpha-n}$ to obtain an $s_2$ and a similar evaluation of
$|\Im|x_1-y|^{\alpha-n}-\Im|x_2-y|^{\alpha-n}|$.
Conclude from this that for any $y\in A$ 
$$
\big||x_1-y|^{\alpha-n}-|x_2-y|^{\alpha-n}\big|\leq c''(\min(|x_1-y|,|x_2-y|))^{-n}|x_1-x_2|.
$$
Choose $K>0$ such that $B(y_1,K)\cap B(y_2,K)\supset A$
for all $y_1,y_2\in\skO$. Then $K$ is independent of $x_1$ and $x_2$ and
\begin{eqnarray*}
&&\int_A(\min(|x_1-y|,|x_2-y|))^{-n}\,dy\\
&\leq&\int_{B(x_1,K)\setminus B_1}
|x_1-y|^{-n}\,dy +\int_{B(x_2,K)\setminus B_2}|x_2-y|^{-n}\,dy\\
&\leq&2\Omega_n(\log K-\log(\frac{2}{3}|x_1-x_2|))\leq C'(1+|x_1-x_2|^{\epsilon-1}),
\end{eqnarray*}
where $C'$ is independent of $x_1$ and $x_2$. The last evaluation
holds because $\epsilon-1<0$.
\end{proof}

\begin{corollary}
\label{cor:alternative}
Let $f\in C_a^l(\Rn)$. Let $\alpha\in\C$ with $1\leq\Re\alpha<\min(a,n)$ be
given. Then 
$$
\Ia f\in C^{(l+[\Re\alpha]-1)+\langle\epsilon\rangle}(\Rn)
$$
for any $0<\epsilon<1$.
\end{corollary}

\begin{proof}
Let $0<\epsilon<1$ be given. If $\Re\alpha=1$ the claim is the
previous lemma. If $\Re\alpha>1$ write $\alpha=\beta+1$. From Proposition \ref{prop:sammensaetning}
$$
I^\alpha f=I^\beta(I^1f).
$$
According to the previous lemma
$I^1 f\in C^{l+\langle\epsilon\rangle}(\Rn)$, so the claim
follows from Corollary \ref{cor:I_decrease}.
\end{proof} 

\begin{proof}[Proof of Theorem \ref{thm:alternativt}]
Use Remark \ref{rem:no_Holder_x} to modify the conclusions of
Proposition \ref{prop:r_alpha} and \ref{prop:I_alpha} regarding the set of
definition of the extension when the H\"older continuity on the
derivatives of $f$ is replaced by ordinary continuity. Use this in
following the lines of the proof of Proposition
\ref{prop:sammensaetning_-k_k}: Let $x\in\Rn$ be given. The map
$$
\alpha\mapsto(I^{\alpha+k}f)(x)
$$
is holomorphic in $\{\alpha\in\C|-k<\Re\alpha<\delta\}$. Since Lemma
\ref{lemma:decay} and Corollary \ref{cor:alternative} ensures, that
$$
I^kf\in C^{(k-1)+\langle\epsilon\rangle}_\delta(\Rn)
$$
for all $0<\epsilon<1$, the map
$$
\alpha\mapsto(I^\alpha(I^kf))(x)
$$
is holomorphic in
$$
\bigcup_{0<\epsilon<1}\{\alpha\in\C|-(k-1+\epsilon)<\Re\alpha<\delta\}=\{\alpha\in\C|-k<\Re\alpha<\delta\}.
$$
Thus by Proposition \ref{prop:sammensaetning} and analytic extension 
\begin{equation} \label{eq:fra_lemma}
I^\alpha(I^k f)(x)=I^{\alpha+k}f(x),
\end{equation}
when $-k<\Re\alpha<\delta$.
The last step of the proof of Proposition \ref{prop:sammensaetning_-k_k}
requires Lemma \ref{lemma:If} the conclusion of which does not hold
for an arbitrary $f\in C_a(\Rn)$
($I^0(f)$ does not necessarily exist). But we can use Lemma
\ref{lemma:alternativt} to replace Lemma \ref{lemma:If} by (see (\ref{eq:If}))
\begin{displaymath}
\lim_{s\to0_+}(I^s f)(x)=
\lim_{s\to0_+}\frac{\Gamma(s)}{H_n(s)}
r^{s-n}(\tau_x f)
=\lim_{s\to-1_+}x^s_+(M_{\tau_x f})=M_{\tau_x f}(0)=f(x).
\end{displaymath}
Thus, using (\ref{eq:fra_lemma}), we have that
\begin{displaymath}
\lim_{s\to-k_+}I^s(I^k f(x))=f.
\end{displaymath} 
This in connection with (\ref{eq:udregning}) proves the theorem.
\end{proof}

\section{The Inversion Formula in Terms of the Laplacian}

It is known, cf. \cite{Helgason}, that if $k$ is even, the inversion formula can be stated by means of the Laplacian,
$\Delta$, instead of the more complicated Riesz potentials. In fact

\begin{theorem}
\label{thm:Laplacian}
When $k$ is even, and
$f\in C^2(\Rn)$, and $f$ and all its first and second order derivatives
are $O(|x|^{-k-\epsilon})$ for some $\epsilon>0$, then
\begin{equation}
\label{eq:Laplace}
  f=(4\pi)^{-\frac{k}{2}}\frac{\Gamma\left(\frac{n-k}{2}\right)}
  {\Gamma\left(\frac{n}{2}\right)}(-\Delta)^{\frac{k}{2}}(\hat{f})\,\check{}.
\end{equation}
\end{theorem}

\begin{proof}
Follow the lines of \cite[p.16-17]{Helgason}: First
notice that it suffices for $f$ to be continuous and
$O(|x|^{-k-\epsilon})$ for some $\epsilon>0$ in order to have formula
\cite[(34)]{Helgason} for the $k$-plane transform; that is,
\begin{equation}
\label{eq:iteration}
(\hat{f})\,\check{}\,(x) = \Omega_k\int^\infty_0 F(r,x)r^{k-1}\,dr,
\end{equation} 
for any $x\in\Rn$, where $F(r,x) = \frac{1}{\Omega_n}\int_{S^{n-1}}f(x+r\omega)\,d\omega$. Here $d\omega$ is the Haar measure on the unit sphere $S^{n-1}$ in
$\R^n$ with total mass 
$\Omega_n=\frac{2\pi^{\frac{n}{2}}}{\Gamma(\frac{n}{2})}$.
Then notice, that the demands on the decay of the derivatives of $f$
allows us to apply the Laplacian (with respect to $x$) on
(\ref{eq:iteration}) by interchanging it with the integration. By
means of Darboux's equation, it can now be seen, as in
\cite[p.16-17]{Helgason}, that
\begin{displaymath}
\Delta((\hat{f})\,\check{}\,)(x)=\left\{
\begin{array}{ll}
-\Omega_k(n-k)f(x) & k=2\\
-\Omega_k(n-k)(k-2)\int^\infty_0F(r,x)r^{k-3}\,dr & k\neq2
\end{array}
\right..
\end{displaymath}
When $k=2$ this is (\ref{eq:Laplace}). For $k\neq2$ the expression is
similar to (\ref{eq:iteration}) - the power of $r$ in the integral has
just been reduced and it is still larger than -1. Thus the Laplacian
can be applied once more without inducing further demands on f or its
derivatives. Continued iteration proves (\ref{eq:Laplace}).
\end{proof}

Can the Theorem \ref{thm:omvendingssaetningen} be used to enlarge the class of functions for which (\ref{eq:Laplace}) holds? Not much, I think. Some relevant thoughts are the following: Let $\alpha_0\in\C$ be given. Using Definition \ref{def:Riesz} and Green's formula it is not hard to see, that for $\phi\in C^2(\Rn)$ with sufficient decay of $\phi$ and all it's first and second order derivatives ($O(|x|^{-2-\epsilon})$ for some $\epsilon>0$ is enough),
\begin{equation}
\label{eq:push}
I^\alpha\Delta\phi = -I^{\alpha-2}\phi
\end{equation}
in some strip $\{\alpha\in\C|2<\Re\alpha<2+\delta\}$. If furthermore $\phi\in C^{l+}(\Rn)$ for some integer $l\geq-\Re\alpha_0+2$, Proposition \ref{prop:I_alpha} can be used to extend both sides of (\ref{eq:push}) holomorphically to $\alpha_0$ and thus prove (\ref{eq:push}) for $\alpha=\alpha_0$.

Iterating (\ref{eq:push}) and then using Lemma \ref{lemma:If} proves that when $k$ is even and positive, and $h\in C^{k+}(\Rn)$, and $h$ and all it's derivatives of order less than or equal to $k$ have a certain decay ($O(|x|^{-2-\epsilon}$) for some $\epsilon>0$ is enough), then
$$
(-\Delta)^{\frac{k}{2}}h = I^{-k}h.
$$
Thus we see from Theorem \ref{thm:omvendingssaetningen} that
(\ref{eq:Laplace}) holds also for $f\in C(k,n)$ when, in stead of
decay demands on derivatives of $f$, we demand a certain decay of $(\hat{f})\,\check{}$ and all it's derivatives of order less than or
equal to $k$ ($O(|x|^{-2-\epsilon})$ is enough). Notice, that since
$(\hat{f})\,\check{}$ is proportional to $I^k f$ the derivatives of
$(\hat{f})\,\check{}$ do exist according to Corollary
\ref{cor:I_decrease}.

\section*{Acknowledgments}

My sincere thanks to Sigurdur Helgason and Boris Rubin for their highly useful
suggestions and comments during my work with this paper. Also thanks to Henrik
Schlichtkrull for all the time he spent commenting on the paper.

\end{document}